\documentclass[a4paper]{article}

\usepackage{amssymb}

\begin{document}

\title{Kronecker's and Newton's approaches to solving~: A first
comparison$^{1}$ \\
{(\Large Extended Abstract)}} 

\author{D. Castro$^{2}$, K. H\"agele$^{3}$, J. E.
Morais$^{4}$, L. M. Pardo$^{2}$}

\addtocounter{footnote}{1}\footnotetext{\small : Research was partially
supported by the Spanish grant~: PB96--0671--C02--02.}

\addtocounter{footnote}{1}\footnotetext{\small : Departamento de
Matem\'aticas, Estad\'{\i}stica y Computaci\'on.\\
{\hbox{\hspace{.7cm}}} Facultad de Ciencias, Universidad de Cantabria,
E-39071 Santander, Spain}

\addtocounter{footnote}{1}\footnotetext{\small : Computer Science 
Department, Trinity College, Dublin 2, Ireland}

\addtocounter{footnote}{1}\footnotetext{\small : Departamento de
Matem\'atica e Inform\'atica, Campus de Arrosad\'{\i}a.\\
{\hbox{\hspace{.7cm}}} Universidad P\'ublica de Navarra, E-31006 Pamplona,
Spain}

\date{}

\maketitle

\smallskip

\medskip
{\em 1991 Math. Subj. Class.:  11D72,11J17,11J61,11J68, 11Y05, 68Q25}



%



%

\typeout{MathMacro: 2/96, update 1/97 kh,lm}


\typeout{Pageformat, }

\setlength{\textheight}{22cm} 

\setlength{\textwidth}{15cm}

\setlength{\topmargin}{0cm}

\setlength{\baselineskip}{10pt}

\setlength{\normalbaselineskip}{10pt} 

\setlength{\parindent}{0mm}

\normalbaselines


\typeout{Spacing, }

\newcommand{\spar} {\vskip 0.25cm}

\newcommand{\mpar} {\vskip 0.50cm}

\newcommand{\bpar} {\vskip 0.75cm}


\typeout{Abbreviations, }

\newcommand{\slp} {straight--line program}

\newcommand{\slps} {straight--line programs}


\typeout{Fonts, }

\font\fiverm=cmr5

\font\sevenrm=cmr7

\font\ninerm=cmr9

\font\tenrm=cmr10

\font\nineit=cmti9

\font\ninebf=cmbx9

\font\ninesl=cmsl9

\font\tensc=cmcsc10


\typeout{Environments, }

\newtheorem{proposition}{Proposition}

\newtheorem{definition}[proposition]{Definition}

\newtheorem{notation}[proposition]{Notation}

\newtheorem{theorem}[proposition]{Theorem}

\newtheorem{proof}[proposition]{Proof} 

\newtheorem{corollary}[proposition]{Corollary} 

\newtheorem{lemma}[proposition]{Lemma}

\newtheorem{remark}[proposition]{Remark} 

\newtheorem{observation}[proposition]{Observation} 

\newtheorem{main}[proposition]{Main Theorem} 

\newtheorem{problem}{Problem}

\newtheorem{example}{Example}

\newcommand{\pref}[1] {{\rm (\ref{#1})}}

\renewcommand{\theenumi}{\roman{enumi}}

\renewcommand{\labelenumi}{\theenumi)}


\typeout{Maths, }


\def\ifm#1#2{\relax\ifmmode#1\else#2\fi}



\newcommand{\klk}    {\ifm {,\ldots,} {$,\ldots,$}}


\newcommand{\plp}    {\ifm {+\cdots+} {$+\ldots+$}}


\newcommand{\absolut}[1]{\ifm {\vert\,{#1}\,\vert} {$\vert\,{#1}\,\vert$}} 

\newcommand{\absolutd}  {\ifm {\absolut{\cdot}} {$\absolut{\cdot}\vert$}} 

\newcommand{\absv}[1]	{\ifm {{\absolut{#1}}_{\nu}} {$\absolut{#1}_{\nu}$}} 

\newcommand{\absvd}	{\ifm {{\absolutd}_{\nu}} {${\absolutd}_{\nu}$}}


\newcommand{\norm}[1]	{\ifm {\vert\vert{#1}\vert\vert} {$\vert\vert{#1}\vert\vert$}}


\newcommand{\undertext}[1] {\ifm {\underline{\hbox{#1}}} {$\underline{\hbox{#1}}$}}


\newcommand{\vect}[2] {\left(\begin{array}{c}{#1}\\

\vdots\\{#2}\end{array}\right)}


\newcommand{\prf}{\spar\noindent {\em Proof}.-- \ } 

\newcommand{\cqfd}{\hfill\vbox{\hrule height 5pt width 5pt }\bigskip}


\typeout{Sets, }


\def \N{{\rm I\kern -2.1pt N\hskip 1pt}} 

\def \Z{{\ifm{{\rm Z\!\!Z}}{${\rm Z\!\!Z}$}}}

\def \R{{\rm I\kern -2.2pt R\hskip 1pt}} 

\def \K {{\rm I\kern -2.1pt K\hskip 1pt}} 

\def \F {{\rm I\kern -2.1pt F\hskip 1pt}} 

\def \P{{\rm I\kern -2.2pt P\hskip 1pt}} 

\def \Fp {\F_p}

\def \wzero {\setminus\{0\}}

\newcommand{\Barq}  {\vrule height0.65em width0.05em depth0em \,}

\newcommand{\Q}        {\ifm {\mbox{\rm Q}\hspace{-0.54em}\Barq\>\;}

                             {$\mbox{\rm Q}\hspace{-0.54em}\Barq\>\;$}}

\newcommand{\C}        {\ifm {\mbox{\rm C}\hspace{-0.45em}\Barq\>}

                             {$\mbox{\rm C}\hspace{-0.45em}\Barq\>$}}

\newcommand{\miniQ}        {{\ninerm Q}

\hspace{-0.40em}\vrule height0.27em width0.025em depth0em \,\>\;}


\typeout{more Maths.}



\newcommand{\om}[2]   {{#1}_1 \klk {#1}_{#2}}


\newcommand{\xo}[1]  {\ifm {\om X {#1}} {$\om X {#1}$}}

\newcommand{\xon}    {\ifm {\om X n} {$\om X n$}}

\newcommand{\yo}[1]  {\ifm {\om Y {#1}} {$\om Y {#1}$}}

\newcommand{\yon}    {\ifm {\om Y n} {$\om Y n$}}

\newcommand{\tto}[1]  {\ifm {\om T {#1}} {$\om T {#1}$}} 

\newcommand{\ton}    {\ifm {\om T n} {$\om T n$}}

\newcommand{\fo}[1]  {\ifm {\om f {#1}} {$\om f {#1}$}}

\newcommand{\fon}    {\ifm {\om f n} {$\om f n$}}

\newcommand{\go}[1]  {\ifm {\om g {#1}} {$\om g {#1}$}}

\newcommand{\gon}    {\ifm {\om g n} {$\om g n$}}

\newcommand{\ho}[1]  {\ifm {\om h {#1}} {$\om h {#1}$}}

\newcommand{\hon}    {\ifm {\om h n} {$\om h n$}}

\newcommand{\iho}[1]  {\ifm {(\om h {#1})} {$(\om h {#1})$}}

\newcommand{\ihon}    {\ifm {\iho n} {$\iho n $}}


\newcommand{\ifo}[1]  {\ifm {(\om f {#1})} {$(\om f {#1})$}}

\newcommand{\ifon}    {\ifm {\ifo n} {$\ifo n $}}

\newcommand{\igo}[1]  {\ifm {(\om g {#1})} {$(\om g {#1})$}}

\newcommand{\igon}    {\ifm {\igo n} {$\igo n $}}

\newcommand{\rifo}[1]  {\ifm {\sqrt{\ifo {#1} }} {$\sqrt{\ifo {#1}}$}}

\newcommand{\rifon}    {\ifm {\rifo n } {$\rifo n $}}


\newcommand{\Ring}[2]  {\ifm {{#1}[#2]} {${#1}[#2]$}}


\newcommand{\Axon}   {\Ring {A} {\xon}}

\newcommand{\Axo}[1]   {\Ring {A} {\xo {#1}}}

\newcommand{\Ayon}   {\Ring {A} {\yon}}

\newcommand{\Ayo}[1]   {\Ring {A} {\yo {#1}}}

\newcommand{\Bton}   {\Ring {B} {\ton}}

\newcommand{\Bxon}   {\Ring {B} {\xon}}

\newcommand{\Byon}   {\Ring {B} {\yon}}

\newcommand{\Kton}   {\Ring {K} {\ton}}

\newcommand{\Kxon}   {\Ring {K} {\xon}}

\newcommand{\Kxo}[1]   {\Ring {K} {\xo {#1}}}

\newcommand{\Kyon}   {\Ring {K} {\yon}}

\newcommand{\Kyo}[1]   {\Ring {K} {\yo {#1}}}

\newcommand{\Rton}   {\Ring {R} {\ton}}

\newcommand{\Rxon}   {\Ring {R} {\xon}}

\newcommand{\Rxo}[1]   {\Ring {R} {\xo {#1}}}

\newcommand{\Ryon}   {\Ring {R} {\yon}}

\newcommand{\Ryo}[1]   {\Ring {R} {\yo {#1}}}

\newcommand{\kton}   {\Ring {k} {\ton}}

\newcommand{\kxon}   {\Ring {k} {\xon}}

\newcommand{\kyon}   {\Ring {k} {\yon}}

\newcommand{\Zton}   {\Ring {\Z} {\ton}}

\newcommand{\Zxon}   {\Ring {\Z} {\xon}}

\newcommand{\Zyon}   {\Ring {\Z} {\yon}}

\newcommand{\Qton}   {\Ring {\Q} {\ton}}

\newcommand{\Qxon}   {\Ring {\Q} {\xon}}

\newcommand{\Qxo}[1]   {\Ring {\Q} {\xo {#1}}}

\newcommand{\Qyon}   {\Ring {\Q} {\yon}}


\newcommand{\bino}[2] {\ifm {{#1}\choose{#2}} {${#1}\choose{#2}$}}


\newcommand{\Tr}  {\mbox{\em Tr}}

\newcommand{\Trt} {\ifm {\widetilde{\mbox{\em Tr}}} {$\widetilde{\mbox{\em Tr}}$} }

\newcommand{\Trb} {\mbox{\rule[.29cm]{0.26cm}{0.1mm}\hspace{-.34cm}\em Tr}}



\section{Introduction}
\label{section-introduction}

Let $K$ be a number field containing the field of Gaussian rationals $\Bbb{Q}[i]\subseteq K$.
In these pages we are mainly interested in the computation of $K-$rational points
of zero--dimensional algebraic varieties given by systems of multivariate 
polynomial equations. Namely, let $f_1,\ldots,f_s\in\Bbb{Z}[X_1,\ldots,X_n]$ be 
a sequence of multivariate polynomials with integer coefficients.
Let $V(f_1,\ldots,f_s)\subseteq \Bbb{C}^n$ be the complex algebraic variety 
of their common zeros, i.e.
$$V(f_1,\ldots,f_s):=\{ x\in \Bbb{C}^n \; : \; f_i(x) =0, \; 1\leq i \leq s\}.$$
For sake of simplicity, let us assume that $V(f_1,\ldots,f_s)$ is a finite 
set (i.e. a zero--dimensional algebraic variety). The set of $K-$rational 
points in $V(f_1,\ldots,f_s)$ is the set of common zeros of the system 
$f_1,\ldots,f_s$ whose coordinates lie in $K^n$, namely
$$V_K(f_1,\ldots,f_s):=  \{x\in K^n\; : \; f_1(x)=\cdots=f_s(x)=0\}.$$
The goal of these pages will be to discuss several aspects of procedures 
performing the following task~: Assume that the field $K$ is fixed.
Given the sequence $f_1,\ldots,f_s$, compute all $K-$rational points 
in $V_K(f_1,\ldots,f_s)$ (or eventually all $K-$rational points in $V_K(f_1,\ldots,f_s)$ 
of bounded height). 

\spar

\noindent Note that the assumption on the field $K$ is not very restrictive~:
For every zero $\zeta\in V(f_1,\ldots,f_s)$, there exists a minimal number field
$K(\zeta)$ containing all the coordinates of $\zeta$. The degree of the field
extension $K(\zeta)$ over $\Bbb{Q}$ can also be denoted by $\deg(\zeta)$. 
In the sequel,  the degree $[K\, :\, \Bbb{Q}]$
may be replaced by  $\deg(\zeta)$, and the results will equally hold.

\spar

\noindent For our study, we consider a precomputation task which prepares the
input $ F:=(f_1,\ldots,f_s)$, before we study the desired $K-$rational points. 
Procedures performing this precomputation task are usually called  
{\em multivariate polynomial system solvers} applied to the input $F$. 
The output of such polynomial system solvers is called the {\em solution} 
of the system $F$. 
Observe that all usual notions of solution of $F$
will yield a description of the variety $V(f_1,\ldots,f_s)$ 
(cf. also \cite{CaGiHeMaPa99}).

\spar

\noindent Here, we consider two (conceptually different) notions
which define what a {\em solution of the system $F$} should be~:
coming from different fields, the notions are related to a symbolic/geometric 
and a numerical analysis/diophantine approximation context~: 
{\em Kronecker's geometric solution} and 
{\em Newton's approximate zero solution}.

\spar

Thus, our study includes a comparative study of both approaches with regards
to the basic problem described above.  It must be said that our study is not
intended to be either complete or definitive. It just tries to point out
some similarities and differences between both approaches to solving that
yield some statements and some open questions of interest. In this sense, we
have tried to write down as many comments as possible to clarify (as much as
we can) the relations between both approaches to solving.

\spar

\noindent Moreover, we have tried to put both approaches under the same hypothesis.
This means, that our input system of multivariate polynomials
$F:=(f_1,\ldots,f_s)$ is well--suited for the application of either Kronecker's or
Newton's approach to solving. 
Therefore we will assume the following hypotheses~:

\label{conditions-systems} 

\begin{itemize}
\item[$i)$] The number of equations equals the number of variables (i.e. $s=n$ above)
\item[$ii)$] The variety $V(f_1,\ldots, f_n)$ is zero--dimensional and contains exactly $D$ 
points,
i.e. the degree of $V(f_1,\ldots,f_n)$ (in the sense of \cite{Heintz83}) 
is exactly $D$.
\item[$iii)$] The $K-$rational points in $V_K(f_1,\ldots,f_n)$ are smooth with
 respect to the system $F:=(f_1,\ldots,f_n)$, i.e. for every $\zeta\in
V_K(f_1,\ldots,f_n)$, the jacobian matrix %
$$DF(\zeta):=\left( {{\partial f_i} \over {\partial X_j}} (\zeta) \right)_{1\leq i,j\leq n}$$
is a non--singular matrix $DF(\zeta)\in GL(n,K)$.
\item[$iv)$] The sequence $f_1,\ldots,f_n\in \Bbb{Z}[X_1,\ldots,f_n]$ is a reduced
regular  sequence, i.e. for every $i$, $1\leq i \leq n-1$, the ideals
$(f_1,\ldots,f_i)$ are radical ideals of codimension $i$ in
$\Bbb{Q}[X_1,\ldots,X_n]$. \item[$v)$] The degrees of the input polynomials satisfy
$\deg(f_i)\leq 2$, for $1\leq i \leq n$.  \end{itemize}

\noindent It must be said that constraints $i)$ and $iv)$ are not relevant
for Kronecker's approach to solving. Applying the iterative version of
Bertini's Theorem (as described in \cite{Morais97,Haegele98,HaMoPaSo99} or
\cite{GiSc99}) we  can easily reduce the over--determined
input system to a system satisfying properties $i)$ and $iv)$. Anyway, we
prefer to keep this hypothesis to simplify exposition, notations, and --
hopefully -- reading.

\spar

\noindent The rest of the introduction presents the new results,
classified into three main categories~:

\begin{enumerate}
\item Newton's approach to solving.

Here, we show how to extend the approximate zero theory introduced by S. Smale
in \cite{Smale81} (cf. also \cite{Smale85,Smale86,Smale86a}), 
and deeply developped in
colaboration with M. Shub in the series of papers 
\cite{ShSm85,ShSm86,ShSm93,ShSm93A,ShSm93B,ShSm94B,ShSm94} and
in~\cite{Malajovich93,Malajovich94,Malajovich95,Yakoubsohn95,Yakoubsohn95b,DeSh1,DeSh2,
Dedieu97,Dedieu97c,Dedieu97b,Dedieu96}, 
to a diophantine approximation context. 

\item Kronecker's approach to solving.

This recalls Kronecker's approach to solving and shows the main statements 
which relate both approaches by means
of an algorithm based on the $L^3$ (or $LLL$) reduction procedure 
(as introduced in \cite{LeLeLo82b} and used in \cite{KaLeLo84,Lenstra84A}). 

\item Application: Computation of splitting field and Lagrange resolvent.

Finally we exhibit an algorithm that combines both approaches to compute
efficiently the splitting field of an univariate polynomial equation and
also the corresponding Lagrange resolvent.

\end{enumerate}
\section{Newton's approach to solving}
\label{subsection-introduction-newton}

\noindent Let $M_K$ be a proper class of absolute values on the number 
field $K$ in the sense of \cite{Lang83}. 
For every $\nu \in M_K$ we have an absolute value
$\vert \cdot \vert_{\nu} ~: K\longrightarrow \mathbb{R}$.
The class $M_K$ is chosen such that it satifies Weil's product formula 
with respect to well-defined multiplicities. 
We denote by $S\subseteq M_K$ the set of sub--indices $\nu\in M_K$ 
such that the absolute value $\absvd$ is archimedean and, consequently, 
by $M_K\setminus S$ the class of sub--indices
$\nu \in M_K$ such that $\absvd$ is non--archimedean.
For every $\nu\in M_K$, we shall denote by $K_{\nu}$ the completion of $K$ 
with respect to the absolute value $\absvd$. We also denote by
$\vert \cdot \vert_{\nu} ~: K_{\nu}\longrightarrow \mathbb{R}$ the corresponding 
extension of $\absvd$ to the completion $K_{\nu}$.

\spar

\noindent Let $\zeta\in V_K(f_1,\ldots,f_n)$ be a smooth $K-$rational 
point of the zero--dimensional complex algebraic variety $V(f_1,\ldots, f_n)$.
We are  interested in approximating $\zeta$ using iterations of the 
Newton operator. 
Therefore, we introduce the Newton operator of system $F$ 
as the following list of rational mappings~:
$$N_F(X_1,\ldots,X_n):= \pmatrix{X_1 \cr \vdots \cr X_n} - 
Df(X_1,\ldots,X_n)^{-1}\pmatrix{f_1(X_1,\ldots,X_n)\cr \vdots \cr
f_n(X_1,\ldots,X_n)}.$$ %
An {\em approximate zero} $z$ in $K^n$ for the system $F$ with associate zero 
$\zeta\in V_K(f_1,\ldots,f_n)$ with respect to the absolute value $\absvd$
is a point  such that the sequence of iterates of the Newton operator is well--defined 
and converges quadratically to $\zeta$. Roughly speaking, an approximate zero 
$z\in K^n$ with associate zero $\zeta\in K^n$ is a point which lies in the basin 
of attraction of the actual zero $\zeta$ with respect to the Newton operator 
$N_F$. Formally, we define approximate zeros as follows~:

\begin{definition}
Let $F:=(f_1,\ldots,f_n)$ be a system of multivariate polynomials with integer coefficients~: 
$f_i\in\Bbb{Z}[X_1,\ldots,X_n]$ for $1\leq i \leq n$. Let
$\nu\in M_K$ define an absolute value $\absvd ~: K\longrightarrow \mathbb{R}$.
Let $\zeta\in V_K(f_1,\ldots,f_n)$ be a smooth $K-$rational point 
(i.e. $DF(\zeta)\in GL(n,K)$).
Let $z:=(z_1,\ldots,z_n)\in K^n$ be an affine point.
We say that $z$ is an approximate zero of the system $F$ with associate zero
$\zeta\in K^n$ with respect to the absolute value $\absvd$, 
if the following properties hold~:
\begin{itemize}
\item $DF(z)\in GL(n,K)$ is a non--singular matrix.
\item The following sequence is well--defined~:
$$z_1:= N_F(z)\in K^n\mbox{, and }z_k:=N_F(z_{k-1})\mbox{ for }k\leq 2.$$
\item For every $k\in \mathbb{N}$, $k\geq 1$, the following inequality holds~:
$$\| z_k -\zeta \|_{\nu} \leq {{1} \over {2^{2^{k-1}}}}\|z -\zeta\|_{\nu},$$
where $\| \cdot \|_{\nu} ~: K_{\nu} \longrightarrow \mathbb{R}$ is the corresponding 
norm associated to the absolute value $\absvd$. 
\end{itemize}
\end{definition} 

\noindent From a computational point of view, 
we want to compute approximate zeros of smooth $K-$rational points and 
we want to write them over a finite alphabet. In particular for every 
smooth $K-$rational zero $\zeta\in V_K(f_1,\ldots,f_n)$  and every absolute
value $\nu\in M_K$, we consider a subfield $L$ of $K$, such that the completion
$L_{\nu}$ of $L$ with respect to the absolute value $\absvd$ contains the entries of $\zeta$, 
namely $\zeta\in L_{\nu}^n$. Thus, we look for approximate zeros $z\in L^n$
with associate zero $\zeta\in L_{\nu}^n$. Let us observe that if the absolute
value $\absvd$ is archimedean, we may fix $L$ to be $L:=\Bbb{Q}[i]$. Moreover,
we are interested in the heights of approximate zeros $z\in L^n$ with actual
zeros $\zeta\in L_{\nu}^n$. In the case where $L= \Bbb{Q}[i]$, 
the height of a point
$z\in \Bbb{Q}[i]^n$ essentially equals its bit length (i.e. the number of tape
cells in a Turing machine required to write down the number list of digits
describing $z$).   In the sequel, we shall therefore identify the
logarithmic height $ht(z)$ and its bit  length. \spar

\noindent A first relevant task consists in stating 
conditions which are sufficient for verifying the property of being 
an approximate zero. 
This is achieved by means of a local condition 
based on a quantity (called $\gamma$), which is essentially yielded
by the Lipschitz constant appearing in the inverse mapping Theorem (cf.
\cite{Demazure89}, Ch. 1, for instance).  These ideas were introduced by S.
Smale in the early eighties  (cf. \cite{Smale81}) and deeply developped in the
series of papers written  by M. Shub and S. Smale \cite{ShSm85} to
\cite{ShSm94}) below). 
\spar

With the same notations as above, let $\nu\in M_K$ be an absolute value on the 
field $K$. We define the {\em quatity } $\gamma$~:
$$\gamma_{\nu}(F,\zeta):=\sup_{k\ge 2}\left\|{(DF(\zeta))^{-1}(D^{(k)}F(\zeta)) \over k!}  
\right\|_{\nu}^{1 \over k-1},$$
where the norm is considered as the norm with respect to the absolute value 
$\absvd$ of the multilinear operator
$$DF(\zeta)^{-1}D^{(k)}F(\zeta)~: \left(K_{\nu}^n\right)^{k}\longrightarrow
K_{\nu}^n.$$
\noindent This quantity yields a locally sufficient condition for 
having an approximate zero. 
This statement is known as the $\gamma-$Theorem and it holds equally true for 
archimedean and non--archimedean absolute values.

\begin{theorem}[$\gamma-$Theorem]
\label{theorem-gamma}
With the same notations and assumptions as before, 
let $F:=(f_1,\ldots,f_n)$ be a sequence of multivariate polynomials with 
coefficients in $K$. Let $\zeta\in V_K(f_1,\ldots,f_n)$ be a smooth 
$K-$rational zero (i.e. $DF(\zeta)\in GL(n,K)$ is a non--singular matrix). 
Let $\absvd~: K\longrightarrow \mathbb{R}_{+}$ be an absolute value on $K$. 
For every $z\in K^n$ 
satisfying the inequality~:
$$\|\zeta-z\|_\nu\gamma_\nu(F,\zeta)\le {3-{\sqrt 7} \over 2 }$$
holds~: $z$ is an approximate zero of the system $F$ with associate zero $\zeta$ 
with respect to the absolute value $\absvd$. 
\end{theorem}  

\noindent The proof of this statement follows step by step the proof of the usual 
$\gamma-$Theorems (cf. the compiled version in \cite{BlCuShSm98}).

\spar

\noindent To establish upper and lower bounds for the bit length of approximate 
zeros, we have established several technical statements. 
One of them is an extension to the non--archimedean case of the 
well--known Eckardt \& Young Theorem \cite{EcYo36}~: 

\spar

\noindent Let $\nu\in M_K$ be an absolute value over $K$ and $K_{\nu}$ 
the completion of $K$ with respect to the absolute value $\absvd$. 
Let us denote by $\Sigma_{\nu}\subseteq {\mathcal M}_n(K_{\nu})$ the variety 
of singular $n\times n$ matrices with entries in $K_{\nu}$. Similarly, let 
$\Sigma$ be the subset of $\Sigma_{\nu}$ of all singular $n\times n$ matrices 
with entries in $K$. Finally, let 
$$d_{\nu}^{(F)}~: {\mathcal M}_n(K_{\nu}) \times {\mathcal M}_n(K_{\nu}) \longrightarrow 
\mathbb{R}_{+}$$
be the Frobenius (also called Hilbert--Weil) metric on ${\mathcal M}_n(K_{\nu})$ 
with respect to the absolute value $\absvd$.  Then, the following Theorem
holds~:

\begin{theorem}[Eckardt \& Young]\label{theorem-eckardt-young}
Let $\nu\in M_K$ be an absolute value. For every non--singular $n\times n$ 
matrix $A\in GL(n,K)$, the following equality holds~:
$$ d_{\nu}^{(F)}(A,\Sigma) =d_{\nu}^{(F)}(A,\Sigma_{\nu})=\inf\{ d_{\nu}^{(F)}(A,M)\;
: \; M \in \Sigma\} = {{1} \over {\| A^{-1} \|_{\nu}}}.$$
\end{theorem}

For every multivariate polynomial $f\in \Bbb{Z}[X_1,\ldots,X_n]$ with integer
coefficients, we define its logarithmic height $ht(f)$ as the logarithm of
the maximum of the absolute values of its coefficients. This notion
introduced, we have the following statement which shows lower bounds for the
bit length of approximate zeros.

\begin{theorem}[Lower Bounds]
\label{theorem-lower-bounds}
Let $f_1,\ldots,f_n\in\Bbb{Z}[X_1,\ldots,X_n]$ be a
sequence of multivariate  polynomials. Let us assume that the following
properties hold~: \begin{itemize}
\item[$i)$] $\max\{ \deg(f_i)\; :\;1\leq i\leq n\}= 2$ ,
\item[$ii)$] $ht(f_i)\leq h$ for $1\leq i\leq n$.
\end{itemize}

Let $\zeta\in V_K(f_1,\ldots,f_n)$ be a smooth $K-$rational point of the system 
$F:=(f_1,\ldots,f_n)$. 
Let $\absvd~:K\longrightarrow \mathbb{R}_{+}$ be an absolute value defined on $K$, and
let $L\subseteq K$ be a number field such that
$\zeta\in L_{\nu}^n$.  Then, for every 
$z\in L^n$, $z\not = \zeta$ satisfying~:
$$\vert\vert z -\zeta
\vert\vert_{\nu} \gamma_{\nu}(F,\zeta) \leq {{3-\sqrt{7}}\over 2}$$
the following inequality holds~:
$$ht(z)\ge {1\over 3[L\, : \, \Bbb{Q}]}\left(\log
\gamma_{\nu}(F,\zeta)-[L\, : \, \Bbb{Q}]\left(5\log n  +
2 h\right)-3\right).$$     
Using Theorem~\ref{theorem-eckardt-young} above, the
following inequality also holds~:
$$ht(z)\ge {1\over 3[L\, :\, \Bbb{Q}]}\left( \log
d_{\nu}^{(F)}(DF(\zeta)^{-1}, \Sigma_{\nu})- [L\; : \;
\Bbb{Q}]\left(7\log n +3h \right) -5\right).$$ 
Moreover, in the case where $L=\Bbb{Q}[i]$ is the field of Gaussian rationals, the
two previous lower bounds may be rewritten as~:
$$ht(z)\ge {1\over 6}\left(\log
\gamma_{\nu}(F,\zeta)-\left(10\log n +4h +3\right) \right)\mbox{, and }$$
$$ht(z)\ge {1\over 6}\left( \log d_{\nu}^{(F)}(DF(\zeta)^{-1},
\Sigma_{\nu})-\left(14\log n +6h + 5\right) 
\right).$$
\end{theorem}

\noindent Let us observe that the ``negative terms" in the previous lower bounds 
are linear in the input length (i.e. the bit length of the input system 
$F:=(f_1,\ldots,f_n)$) whereas the ``positive part" depends semantically 
on the smooth $K-$rational solution $\zeta\in V_K(f_1,\ldots,f_n)$.

\spar
\noindent Last, but not least, we may also show a few lower bounds 
for the average height of approximate zeros  associated to a $\Bbb{Q}-$definable
irreducible component of the solution variety $V(f_1,\ldots,f_n)$. To this
end, we introduce some additional notations.
Let
$f_1,\ldots,f_n\in \Bbb{Z}[X_1,\ldots,X_n]$ be a sequence of multivariate
polynomials satisfying the hypotheses $i)$ to $v)$ above. 
Let $\zeta\in
K^n$ be a smooth $K-$rational zero of the system $F:=(f_1,\ldots,f_n)$.
Let $V:=V(f_1,\ldots,f_n)\subset \Bbb{C}^n$ be the algebraic variety given as
the common zeros of the polynomials $f_1,\ldots,f_n$. Let
$\mathcal{V}_{\zeta}\subseteq V$ be  the $\Bbb{Q}-$definable irreducible
component of $V$ that contains $\zeta$ .  Let us assume
$D:=deg(\mathcal{V}_{\zeta})$ be the number of points in $\mathcal{V}_{\zeta}$.
Let us observe that $D=deg(\zeta) \leq [K\; : \Bbb{Q}]$. Let us assume 
$$V_{\zeta}:=\{\zeta_1,\ldots,\zeta_D\}.$$
Let $\| \cdot \|~: K^n \longrightarrow \mathbb{R}$ be the standard hermitian norm
induced in $K^n$ by the inclusion $i ~: K\hookrightarrow \Bbb{C}$. A sequence of
points $z:=(z_1,\ldots,z_D) \in \Bbb{Q}[i]^{nD}$ is said to be an approximate zero
of the system $F$ with associate variety $\mathcal{V}_{\zeta}$ that satisfies
the $\gamma-$Theorem if for every $i$, $1\leq i \leq D$, the following holds~:
$$\| z_i - \zeta_i \| \leq {{3-\sqrt{7}}\over {2\gamma(F,\zeta_i)}},$$
where $\gamma(F,\zeta_i)$ is the quantity associate to the hermitian
norm $\| \cdot \|$.

\noindent For every given approximate zero $z:=(z_1,\ldots,z_D)\in \Bbb{Q}[i]^{nD}$
of the system $F$ with associate variety $V_{\zeta}$, the average height (also
the average bit length) of $z$  is defined in the following terms
$$ht_{av}(z):= {{1}\over {D}} \sum_{i=1}^D ht(z_i).$$
Finally, let us denote by $\Bbb{Z}_{K}\subset K$ the ring of algebraic integers
of the number field $K$.
Then, we have the following lower bound for the average bit length of
approximate zeros with associate variety $\mathcal{V}_{\zeta}$~:

\begin{proposition}
\label{proposition-integer-lower}
With the previous notations, let $\zeta\in V_K(f_1,\ldots,f_n)$ be a smooth
$K-$rational with entries in $\Bbb{Z}_{K}$, i.e. $\zeta\in \Bbb{Z}_K^n$.  Let us
also assume that for every archimedean absolute value $\absvd$ (i. e. $\nu \in
S$), the following holds~:
$$3 \|\zeta\|_{\nu} \gamma_{\nu}(F,\zeta)\geq
3-\sqrt{7}.$$ 
Then the average height of any approximate zero $z\in
\Bbb{Q}[i]^{nD}$ of the system $F$ with associate variety  $V_{\zeta}$
satisfies the following inequality~:
$$ht_{av}(z) \geq {1 \over 2} \left[ ht(\zeta)-( {1 \over 2} \log n + \log
2)\right].$$
\end{proposition}

\spar

In order to illustrate the meaning of this lower bound, we give here a few
Corollaries.

\begin{corollary}
\label{corollary-integer-gamma}
With the same notations as in
Proposition~\ref{proposition-integer-lower} above, let 
$\zeta\in \Bbb{Z}_K\cap V_K(f_1,\ldots,f_n)$ be a smooth $K-$rational zero of
the system $F:=(f_1,\ldots,f_n)$ and let us assume that for every archimedean
absolute value $\absvd~: K \longrightarrow \mathbb{R}$ (i. e. for every $\nu
\in S$), the following holds~:
$$\gamma_{\nu}(F,\zeta)\geq 3-\sqrt{7}.$$
Then the average height of any approximate zero $z\in
\Bbb{Q}[i]^{nD}$ of the system $F$ with associate variety  $V_{\zeta}$
satisfies the following inequality~:
$$ht_{av}(z) \geq {1 \over 2} \left[ ht(\zeta)-( {1 \over 2} \log n + \log
2)\right].$$
\end{corollary}

\noindent Moreover, the previous techniques show how to deform a given system
of multivariate polynomials by means of a single additional equation 
of low degree in such a way that the  average bit length of the new system is
essentially greater than the height of the actual zero you want to approximate.

\begin{corollary}
\label{corollary-integer-deformation}
Let $F:=(f_1,\ldots,f_n)$ be a system of multivariate polynomials with integer
coefficients satisfying the conditions $i)$ to $v)$ given on page
\pageref{conditions-systems}.
Let $\zeta=(\zeta_1,\ldots,\zeta_n) \in \Bbb{Z}_K\cap V_K(f_1,\ldots,f_n)$ be a
smooth $K-$rational zero whose coordinates are algebraic integers. Let us now
define the system of polynomial equations in $n+1$ variables~:
$$G:=(g_1,\ldots,g_{n+1})\in \left(\Bbb{Z}[X_1,\ldots,X_{n+1}]\right)^{n+1},$$
given by the following rules~:
\begin{itemize}
\item $g_i:= f_i \in \Bbb{Z}[X_1,\ldots,X_{n+1}]$ for every $i$, $1\leq i \leq n$, 
\item $g_{n+1}:= \left(X_{n+1} -X_n\right) \left(X_{n+1}- (X_n +1)\right).$
\end{itemize}  
Let $\zeta'\in V_K(g_1,\ldots,g_{n+1})\cap \Bbb{Z}_K$ be the affine point given
by~:
$$\zeta':=(\zeta_1,\ldots,\zeta_n,\zeta_n)\in \Bbb{Z}_K^{n+1}.$$
Let $\mathcal{V}_{\zeta'}\subseteq V(g_1,\ldots, g_{n+1})$ be the
$\Bbb{Q}-$definable irreducible component of $V(g_1,\ldots,g_{n+1})$
containing $\zeta'$. Then, the average height of any approximate zero $z\in
\Bbb{Q}[i]^{(n+1)D}$ of the system $F$ with associate variety 
$\mathcal{V}_{\zeta'}$ satisfies the following inequality~:
$$ht_{av}(z) \geq {1 \over 2} \left[ ht(\zeta)-( {1 \over 2} \log (n+1) +
\log 2)\right].$$
\end{corollary}
  
\noindent  We can construct 
several examples  where all of the previous lower bounds for the bit length of
approximate zeros apply. For instance we give here an example 
inspired by a classical univariate example due to M.
Mignotte (cf. \cite{Mignotte89})~:

\begin{example} [Using $\log\gamma$ as in Theorem \ref{theorem-lower-bounds}]
Let us consider
the system of multivariate polynomials $F:=(f_1,\ldots,f_{n+1})$ given
by the following rules~:

\begin{itemize}
\item $f_1:= X_1-2$,
\item $f_i:= X_i-X_{i-1}^2$ for every $i$, $2\leq i \leq n-1$,
\item $f_n:= X_{n+1}-X_n^2$,
\item $f_{n+1}:=X_{n+1}X_n -2(X_{n-1}X_n-1)^2$.
\end{itemize}

This system $F$ has three solutions in $\mathbb{C}^{n+1}$, where two of
them, say $\zeta_1, \zeta_2\in \mathbb{R}^{n+1}$, satisfy the following inequality~:
$$\| \zeta_1 - \zeta_2 \| \leq {{2} \over { 2^{ {2^{n-2} 5}\over {2}}}}\leq {
{2} \over { 2^{ 2^{n-1}}}}.$$
Thus, using the well-known relation between the quantity $\gamma$ and
the separation of roots, we may conclude~:
$${{2} \over { 2^{ 2^{n-1}}}}\geq\|\zeta_1 - \zeta_2 \|\geq {{3-\sqrt{7}}\over
{2\gamma(F,\zeta_i)}}.$$
By  Theorem \ref{theorem-lower-bounds}, we
conclude that for all approximate zeros $z_1, z_2\in \mathbb{Q}[i]^{n+1}$
of the system $F$ associated to $\zeta_1, \zeta_2$ respectively and satisfying
the corresponding $\gamma-$Theorem, the following holds~:
$$ht(z_i) \geq  {{1} \over {6}}\left ( 2^{n-1} -2\log (n+1) \right) -
O(1).$$
\end{example}

Therefore this example (an several others wich should be constructed) allows
us to conclude the following Corollaries~\ref{corollary-extime-p} to
~\ref{corollary-continued-fractions}.

\begin{corollary}\label{corollary-extime-p}
Computing approximate zeros in $\Bbb{Q}[i]$ for archimedean absolute values, using
binary encoding of the output requires exponential running time and
exponential output length, and these two lower bounds cannot be improved with
this encoding. Namely,  computing approximate zeros with binary encoding is in
the complexity class {\bf EXTIME} $\setminus$ {\bf P}. 
\end{corollary}

\begin{corollary}\label{corollary-floating}
Floating point encoding of approximate zeros requires exponential
number of digits and this lower bound cannot be improved. Namely,  
floating point encoding does not suffice to compute approximate zeros of
systems of
multivariate polynomial equations.  
\end{corollary}

\noindent As continuous fraction encoding of numbers in $\Bbb{Q}[i]$ is close
to the binary encoding, we easily conclude the following~:

\begin{corollary}\label{corollary-continued-fractions}
Computing approximate zeros in $\Bbb{Q}[i]$ for archimedean absolute values, using
continuous fraction  encoding of the output requires exponential running
time and exponential output length, and these two lower bounds cannot be
improved with this encoding. Namely,  computing approximate zeros with
continuous fraction encoding is in the complexity class {\bf EXTIME}
$\setminus$ {\bf P}. 
\end{corollary}

\noindent These lower bounds suggest 
that a central point of interest should be to study
the bit length of approximate zeros satisfying the $\gamma-$Theorem. 
In order to shed some light in this direction, we prove the following
statements~:

\begin{theorem}[Upper Bounds]
\label{theorem-upper-bounds} 
Let $f_1,\ldots,f_n\in \Bbb{Z}[X_1,\ldots,X_n]$ be polynomials with integer 
coefficients. Let us assume that the following properties hold~:

\begin{itemize}
\item $\max\{\deg(f_i): 1\leq i \leq n\} \leq 2$, and
\item $ht(f_i)\leq h$ for $1\leq i \leq n$.
\end{itemize}

Let $\zeta\in V_K(f_1,\ldots,f_n)$ be a smooth $K-$rational point.
Let $\absvd ~: K \longrightarrow \mathbb{R}_+$ be an absolute value on $K$. 
Then, the following inequality holds~:
$$\log\gamma_{\nu}(F,\zeta)\leq 3[K\; :\;\Bbb{Q}]
n\left(n^2+ 4\log n +  h +ht(\zeta)+3\right).$$
\end{theorem}

\noindent In particular, we show the following estimate for the bit length of
approximate zeros in $\Bbb{Q}[i]^n$~:
\spar


\begin{corollary}[Upper bound on the bit length of approximate zeros]
\label{corollary-upper-length}
With the same assumptions and notations as in Theorem
\ref{theorem-upper-bounds}  above, let $\zeta\in
V_K(f_1,\ldots,f_n)$ be a smooth $K-$rational zero, and let $\absvd$ be an
absolute value on $K$. Let $L\subseteq K$ be a number field such that
$\zeta\in L_{\nu}^n$. Then there exist approximate zeros $z\in L^n$ of the
system $F:=(f_1,\ldots,f_n)$ with approximate zero $\zeta$ with respect to the
absolute value $\absvd$, such that the logarithmic height $ht(z)$ of $z$ is at
most linear in the following quantities~:  %
$${{1}\over {[L\, : \, \Bbb{Q}]}} \log \vert \Delta_L \vert + [K\, : \, \Bbb{Q}]
n\left( n^2 +h +nht(\zeta) \right),$$
where $\vert \Delta_L\vert $ is the absolute value of the discriminant of the
field $L$.
\spar

Moreover, in the case where $L=\Bbb{Q}[i]$ (for
instance, if $\absvd$ is archimedean),
there exist approximate zeros $z\in \Bbb{Q}[i]^n$ for the system $F$ with associate
zero $\zeta$ with respect to $\absvd$ such that their bit length is at most
linear in the following quantity~:
$$[K\, : \, \Bbb{Q}] n\left( n^2 +h +nht(\zeta) \right)\mbox{, in other words~:}$$
$$ht(z)\leq O\left( [K\, : \, \Bbb{Q}]
n\left( n^2 +h +nht(\zeta)\right)\right).$$
\end{corollary}

\noindent Let us observe, that these two upper bounds above 
(i.e. Theorem~\ref{theorem-upper-bounds} and Corollary
\ref{corollary-upper-length})  depend mainly on the input
length $n^2h + n^2$ and on two parameters  which in turn depend on the actual
zero to approximate~:  the degree $[K\, : \, \Bbb{Q}]$ of a field containing
the coordinates of the zero  and the logarithmic height of the zero
$ht(\zeta)$.  These two quantities are bounded respectively by the geometric
B\'ezout inequality  (cf. \cite{Heintz83} or \cite{Fulton84,Vogel84}) and 
the arithmetic B\'ezout inequality (cf.
\cite{BoGiSo93,Philippon91,Philippon94,Philippon95} or 
\cite{KrPa94,KrPa96,Sombra98,Haegele98,HaMoPaSo99}, for instance).  Moreover,
combining these two upper bounds (Theorem~\ref{theorem-upper-bounds} and
Corollary~\ref{corollary-upper-length}) with the previously shown lower
bounds and several examples,  we may conclude that the upper bounds shown in
Theorem~\ref{theorem-upper-bounds} and Corollary
\ref{corollary-upper-length} are optimal.

\spar

\noindent On the other hand, the 
$\gamma-$Theorem above has some aesthetic consequences which we may explain in
terms of the existence of a {\em universal radius of convergence} independent
of the aboslute value under consideration. To this end, we recall the
well--known Implicit Function Theorem for complete noetherian local rings in
the following terms~:

\begin{theorem} [Non--archimedean Basin of Attraction]
\label{theorem-na-basin}
Let $F:=(f_1,\ldots,f_n)\in \Bbb{Z}[X_1,\ldots,X_n]^n$ be a system of multivariate polynomials 
satisfies the hypotheses of Theorem~\ref{theorem-upper-bounds} above.
Let $\nu\in M_K$  define a  non--archimedan absolute value $\absvd$ on $K$. Let us 
also assume that the restriction 
$$\vert \cdot \vert_{\nu} ~: \Bbb{Q} \longrightarrow \mathbb{R}_{+}$$
defines  a $p-$adic absolute value, where $p\in \mathbb{N}$ is a prime number.  Let
$\zeta\in K^n$ be a smooth $K-$rational zero of the system which lies in the closed unit
shere of $K^n$, i.e.
$$\zeta\in B_{\nu}(0,1):=\{ x\in K^n \; : \; \|x \|_{\nu}\leq 1\}.$$
Let us finally assume that $\vert \det DF(\zeta)\vert_{\nu}=1$. Then, for every 
$z\in B_{\nu}(0,1)$ 
satisfying
$$\| z - \zeta \| _{\nu} \leq {1 \over p}$$
holds~: $z$ is an approximate zero of the system $F$ with associate zero $\zeta$
with respect to the absolute value $\absvd$.
\end{theorem}

\noindent This statement is nothing but the usual Hensel Lemma
in local algebra (cf. \cite{ZaSa58,Morais97}, for instance).   
However,  this statement has a drawback~:
The radius of the basin of attraction centered at $\zeta$ depends on the 
concrete absolute value $\absvd$. The
$\gamma$--Theorem~\ref{theorem-upper-bounds}  above shows that
there exists a {\sl universal radius}, which depends only on the  system $F$
and  the smooth $K-$rational zero, but does not depend on the absolute value. 

\spar

In order to prove this claim, 
let us introduce the following quantity 
$\widetilde{\gamma}(F,\zeta)$. With the same notations and assumptions as 
above, let us define the {\em universal quantity}
$$\widetilde{\gamma}(F,\zeta):= 
\left( \prod_{\nu\in M_K} \max \{ 1,
\gamma_{\nu}(F,\zeta)\}^{n_{\nu}}\right)^{1\over [K:\mathbb{Q}]}.$$ 
Let us observe, that this quantity is well--defined and finite 
according to Theorem~\ref{theorem-upper-bounds} above. Moreover, it does not
depend on any particular absolute value under consideration. 
Thus, we may conclude the following Theorem~:

\begin{corollary}[Universal $\gamma-$Theorem]
\label{theorem-universal-gamma}
With the same notations and assumptions as in Theorem
\ref{theorem-lower-bounds},  for every $z\in \Bbb{Q}[i]^n$ and
every absolute value $\absvd$   satisfying the following inequality
$$\| z - \zeta\|_{\nu} \widetilde{\gamma}(F,\zeta) \leq {{3-\sqrt{7}}\over {2}}$$
holds~: $z$ is an approximate zero for the system $F$ with associate zero $\zeta$ 
and with respect to the absolute value $\nu\in M_K$.
\end{corollary}

\noindent Let us point out that the existence of such a universal 
quantity does not imply the existence of a universal basin of attraction
independent of the absolute value under consideration. In fact, we show the
following (expected) statement~:

\begin{corollary}
\label{corollary-intersection}
Let $F:=(f_1,\ldots,f_n)$ be  a sequence of multivariate polynomials with
integer coefficients satisfying conditions $i)$ to $v)$ of page
\pageref{conditions-systems}.
Let $\zeta\in V_K(f_1,\ldots, f_n)$ be a smooth $K-$rational zero.
The only point $z\in K^n$ that satisfies the universal
$\gamma-$Theorem near $\gamma$ for all absolute values in $M_K$ 
is $z=\zeta$. Namely, for every $z\in K^n$ satisfying the following inequality for 
every $\nu\in
 M_k$
$$\| z- \zeta \|_{\nu} \leq {{3-\sqrt{7}}\over {2\widetilde{\gamma}(F,\zeta)}}$$
holds $z= \zeta$.
\end{corollary}

\section{Kronecker's approach to solving}

In \cite{Kronecker82}, Kronecker introduced a notion of {solution of unmixed
complex algebraic varieties}, which we are going to reproduce here. 
Let $f_1,\ldots,f_i\in\Bbb{Z}[X_1,\ldots,X_n]$ be  a sequence of polynomials 
defining a radical ideal $(f_1,\ldots,f_i)$ of codimension $i$ in 
$\Bbb{C}[X_1,\ldots,X_n]$. Let $V:=V(f_1,\ldots,f_i)\subseteq \Bbb{C}^n$ be the complex 
algebraic variety of codimension $i$ given by the common zeros of the $f_i$. 
A {\em solution} of $V$ is a birational isomorphism of $V$ with some complex 
algebraic hypersurface in a space of adecuate dimension. 

\spar

Technically, this is is expressed as follows. First of all, let us assume that
the variables $X_1,\ldots,X_n$ are in Noether position with respect to the 
variety $V$, i.e. we assume that the following is an integral ring extension~:
$$\Bbb{Q}[X_1,\ldots,X_{n-i}]\hookrightarrow \Bbb{Q}[X_1,\ldots,X_n]/(f_1,\ldots,f_i).$$
Let $u:=\lambda_{n-i+1}X_{n_i+1} + \cdots + \lambda_n X_n\in \Bbb{Q}[X_1,\ldots,X_n]$
be a linear form  in the dependent variables $\{X_{n-i+1},\ldots,X_n\}$. 
Thus we have a linear projection
$${\cal U}~:\Bbb{C}^n \longrightarrow \Bbb{C}^{n-i+1}$$
$$(x_1,\ldots,x_n) \longmapsto \left(x_1,\ldots,x_{n-i}, u(x_1,\ldots,x_n)\right).$$
Let us also consider the restriction
$${\cal U}\mid_V~: V \longrightarrow \Bbb{C}^{n-i+1}.$$
The linear form $u$ is called a {\em primitive element} if and only if
the projection ${\cal U}\mid_V$ defined a birational isomorphism of $V$
with some complex hypersurface $H_u$ in $\Bbb{C}^{n-i+1}$ with minimal equation 
$\chi_u\in\Bbb{Q}[X_1,\ldots,X_{n-i},T]$. Then, a Kronecker solution of variety $V$
is a description of  the primitive element $u$, 
the hypersurface $H_u$ through the minimal equation $\chi_u$, 
and a description of the inverse of the birational isomorphism, 
i.e.  $({\cal U}\mid_V)^{-1}$. 
Formally, this list of items may be described as follows~:

\begin{itemize}
\item The list of variables in Noether position $X_1,\ldots,X_n$, which includes a 
description of the 
dimension of $V$.
\item The primitive element  $u:=\lambda_{n-i+1}X_{n_i+1} + \cdots + \lambda_n X_n$ given 
by its 
coefficients in $\Bbb{Z}$.
\item The minimal equation of the hypersurface $H_u$, namely
$$\chi_u\in \Bbb{Z}[X_1,\ldots,X_{n-i},T].$$
\item A description of $({\cal U}\mid_V)^{-1}$. This description is given by the following 
list~:
\begin{itemize}
\item A non--zero polynomial $\rho\in\Bbb{Z}[X_1,\ldots,X_{n-i}]$.
\item A list of polynomials $v_j\in\Bbb{Z}[X_1,\ldots,X_{n-i},T]$, $n-i+1\leq j \leq n$, 
such that the degrees with respect to variable $T$ satisfy 
$\deg_T(v_j)\leq \deg_T(\chi_u)$, for every $j$, $n-i+1\leq j \leq n$.
\end{itemize} 
such that the following holds
$$({\cal U}\mid_V)^{-1}(x,t):=\left(x_1,\ldots,x_{n-i},\rho^{-1}(x)v_{n-i+1}(x,t), \ldots, 
\rho^{-1}(x)v_n(x,t)\right),$$
where $x:=(x_1,\ldots,x_{n-i})\in \Bbb{C}^{n-i}$ and $t\in \Bbb{C}$.

\end{itemize}

\noindent Kronecker conceived an iterative procedure to solve multivariate 
systems of equations $F:=(f_1,\ldots,f_n)$ defining zero--dimensional 
complex varieties, which can be described in the following terms~: 

\spar

First, you start with system $(f_1)$ and you ``solve" the unmixed variety of 
codimension $1$, $V(f_1)\subseteq \Bbb{C}^n$. Then you proceed iteratively~:
From Kronecker's solution of the variety $V(f_1,\ldots,f_i)$ 
you eliminate the new equation $f_{i+1}$ to obtain a Kronecker solution of 
the ``next" variety $V(f_1,\ldots,f_{i+1})$.
Proceed until you reach $i=n$.
This iterative procedure has two main drawbacks which can be explained in the
following terms~:

\begin{itemize}
\item First of all, the space problem arising with the 
representation of the intermediate polynomials.
The polynomials $\chi_u$, $\rho$ and $v_j$ are polynomials of high degree
(eventually of degree $2^i$) involving several variables. 
Thus, to represent them, one has to handle their coefficients, which amounts 
to the following quantities
$${{2^i+n-i+1} \choose {n-i+1}},$$
which for $i:= n/2$ amounts to more than $2^{n^2/4}$ coefficients of great bit length.

\item Secondly, Kronecker's iterative procedure introduced a nesting of interpolation
procedures required for the iterative process and the linear change of
cooerdinates required by every computation of the Noether normalisation.
This nesting of interpolation procedures is difficult to avoid and increases
the running time complexity.
\end{itemize}

Therefore, the procedure was forgotten by contemporary mathematicians
and hardly mentioned in the literature of algebraic geometry.
Macaulay quotes Kronecker's procedure in \cite{Macaulay16} and so does
K\"onig in \cite{Konig03}. But both of them thought, 
that this procedure requires too much running time to be efficient, 
and it was progressively forgotten. 
Traces of this procedure can be found spread over the algebraic geometry 
literature without giving to it the required relevance. 
For example, Kronecker's notion of solution was used by O. Zariski in 
\cite{Zariski95} to define dimension of algebraic varieties, 
claiming that it was also used in the same form by Severi and others.


\spar

\noindent In 1995, two works rediscovered Kronecker's approach to solving
without previous knowledge of it's existing ancestors. These two works
\cite{GiHeMoPa95} and \cite{Pardo95} were able to overcome the first drawback
(space problem of representation) of the previous methods.  The technical
trick was the use of a data structure coming from semi--numerical
modelling~: straight--line programs.  This idea of representing polynomials
by programs that evaluate them goes back to former works of the same
research group (as \cite{GiHe93,FiGiSm93,KrPa94} or \cite{KrPa96}). Moreover,
these ideas were the natural continuation of the ideas previously developped
in \cite{GiHe91}.

\spar

\noindent To overcome the second drawback (Nesting),
the authors introduced a method based on Newton's 
method applied in a non--archimedean context 
(the approximate zeros in the corresponding
non--archimedean basin of attraction were called {\em Lifting Fibers} in 
\cite{GiHaHeMoMoPa97}). This result was obtained in the two papers
\cite{GiHaHeMoMoPa97,GiHeMoMoPa98}.  The key trick to avoid the nesting
of interpolation procedures was based on Hensel's Lemma (also Implicit Mapping
Theorem). Perhaps, the following statement could explain the relations
existing between Hensel's Lemma and Approximate Zero Theory.
To this end, let us introduce a few more notations.
Let $f_1,\ldots,f_r\in \Bbb{C}[X_1,\ldots,X_n]$ be a sequece of polynomials
defining
a radical ideal of codimension $r$ in $\Bbb{C}[X_1,\ldots,X_n]$.
Let us assume that the variables $X_1,\ldots,X_n$ are in Noether position with
respect to the ideal $I:=(f_1,\ldots,f_r)$, i.e. assume that the following
ring extension is integral
$$\Bbb{C}[X_1,\ldots,X_{n-r}]\hookrightarrow
\Bbb{C}[X_1,\ldots,X_n]/I.$$
Let $P:=(p_1,\ldots,p_{n-r})\in \Bbb{C}^{n-r}$ be an affine point, let ${\mathcal
O}_P$ be the ring of formal power series at $P$, and let ${\cal M}_P$ be
the field of fractions of ${\mathcal O}_P$. Then, the following is finite ring
extension 
$${\mathcal M}_P \hookrightarrow B:={\cal M}_P[X_{n-r+1},
\ldots,X_n]/(f_1,\ldots,f_r),$$
and $B$ is a zero--dimensional ${\cal M}_P-$algebra. Thus, it has some
sense to look for approximate zeros of the solutions in ${\cal M}_P^r$ of
the system of polynomial equations $F:=(f_1,\ldots,f_r)$.
The following statement about Hensel's Lemma explains the connections existing
between Kronecker's solving and Approximate Zero Theory.

\begin{theorem} [Hensel's Lemma]
With the same assumptions and notations as above, let $\zeta\in {\mathcal
M}_P^r$ be a solution of the system $F$. Let
$\| \cdot \| ~: {\cal M}_P^r\longrightarrow \mathbb{R}$ be usual non--archimedian norm
in ${\cal M}_P^r$. Let $\Bbb{C}(X_1,\ldots,X_{n-r})$ be the field of rational
functions. Then,  for every $z\in\Bbb{C}(X_1,\ldots,X_{n-r})^r$,  if  
$\|z\|\leq 1$, and 
$$\|z -\zeta \| < {1 \over 2},$$ 
then $z$ is an approximate zero for the system $F:=(f_1,\ldots,f_r)$ with
associate zero $\zeta\in {\cal M}_P^r$. \end{theorem}

Unfortunately, those two works \cite{GiHaHeMoMoPa97,GiHeMoMoPa98}
introduced (for the Lifting Fibers) run time requirements which 
depend on  the heights of the intermediate varieties 
(in the sense of \cite{BoGiSo93,Philippon91,Philippon94,Philippon95,Sombra98}). 
This drawback was finally overcome in the paper 
\cite{GiHeMoPa97}, 
where integer numbers were represented by straight--line programs and 
the following result established~:

\begin{theorem}\cite{GiHeMoPa97}
\label{theorem-kronecker-solution}
There exists a bounded error probability Turing
machine $M$ which performs the following task~: 
Given a system of multivariate polynomial equations $F:=(f_1,\ldots, f_n)$, 
satisfying the following properties
\begin{itemize} 
\item $\deg(f_i)\leq 2$ and $ht(f_i)\leq h$ for $1\leq i \leq n$,
\item the ideals $(f_1,\ldots,f_i)$ are radical ideals of codimension $i$
in the ring $\Bbb{Q}[X_1,\ldots,X_n]$ for $1\leq i \leq n-1$,
\item the variety $V(f_1,\ldots,f_n)\subseteq\Bbb{C}^n$ is a zero--dimensional
complex algebraic variety.
\end{itemize}

Then, the machine $M$ outputs a Kronecker solution of the variety 
$V(f_1,\ldots,f_n)$. The running time of the machine $M$ is polynomial in 
the following quantities
$$\delta(F) n h,$$
where $\delta$ is the maximum of the degrees of the intermediate varieties 
(in the sense of \cite{Heintz83}), namely
$$\delta(F):=\max\{ deg(V(f_1,\ldots,f_i)) \; : 1\leq i \leq n-1\}.$$
\end{theorem} 

\noindent It must be said that the coefficients of the polynomials involved in a
Kronecker solution of the variety $V(f_1,\ldots,f_n)$ are given by
straight--line programs that evaluate integer numbers. 
However, the complexity estimates for the Turing
machine $M$ are independent from the height.

\spar

\noindent Our attempt in these pages is to compare this 
approach to solving developped by Kronecker 
to that of Newton as described in the previous Subsection.

\spar

\noindent The exposition of new results starts with a small improvement of
the Witness Theorem in \cite{HeSchn80} and \cite{BlCuShSm96} (cf. also
\cite{BlCuShSm98}). When dealing with straight--line program data structures,
some relevant technical methods of comparison must be developped. 
These methods are known as probabilistic zero tests for polynomials 
given by straight--line programs. Examples of these tests are those introduced in 
\cite{Schwartz79,Zippel79,HeSchn80} and the Witness Theorem, introduced in
\cite{HeSchn80} 
for the case of polynomials with integer coefficients, and in 
\cite{BlCuShSm96} for polynomials with coefficients in a number field. 

\spar

As we already had to introduce a few technical notions and methods spread over the
literature of number theory, numerical analysis, algebraic complexity theory
and elimination theory, 
we can give without introducing additional material
for free the following improvement of the estimates for the Witness
Theorem.

\begin{theorem}[Witness Theorem]
\label{theorem-witness}
Let $P\in K[X_1,\ldots,X_n]$ be a non--zero polynomial evaluable
by a non--scalar straight--line program $\Gamma$ of size $L$, non--scalar
depth $\ell$ and parameters in ${\cal F}\subseteq K$.
Let $\omega_{0}\in K$ be  such
that the following holds~:
$$ht(\omega_{0})\geq\max\{\log 2, ht({\cal F})\}.$$
Let $N\in\mathbb{N}$ be a non--negative integer such that
$$\log_{2}N>\log_2 (\ell+1)+(\ell+2)\left(\log_2\log_2(4 L)\right).$$
Let us define recursively the following sequence of algebraic numbers (known
as Kronecker's scheme)~:
$$\omega_{1}:=\omega_{0}^{N}\mbox{, and for }2\leq i \leq n, \mbox{ }
\omega_{i}:=\omega_{i-1}^{N}.$$
Then, the point $\underline{\omega}:=(\omega_{1},\ldots,\omega_{n})\in K^{n}$ 
is a witness for $P$ (i.e. $P(\underline{\omega})\not=0$).
\end{theorem}

\noindent Moreover, we observe that {\em approximate zeros are succint encodings
of generic points of the variety} $V(f_1,\ldots,f_n)$. 
This means that for every
smooth $K-$rational zero $\zeta \in V_K(f_1,\ldots,f_n)$,
the binary encoding of an approximate zero $z\in \Bbb{Q}[i]^n$ is sufficient
information to compute
the $\Bbb{Q}-$irreducible component of $V(f_1,\ldots,f_n)$ containing $\zeta$. 
In more precise terms we show the following statement~:

\begin{theorem} [From Approximate Zeros to Geometric Solution]
\label{theorem-newton-kronecker}
With the same assumptions as in Theorem~\ref{theorem-kronecker-solution}
above, there exists a bounded error probability Turing machine $M$, 
such that taking as input the binary encoding of an approximate zero 
$z\in \Bbb{Q}[i]$ of the system $F$ with associate zero $\zeta\in V_K(f_1,\ldots,f_n)$ 
for an archimedean absolute value $\absvd$ (where $\nu\in S$), $M$ 
outputs a Kronecker solution of the $\Bbb{Q}-$irreducible component $W$ of 
$V(f_1,\ldots,f_n)$ containing $\zeta$. 
Moreover, the running time of this probabilistic Turing machine is polynomial 
in the following quantities
$$\deg(W)\left(n\; h\; ht(z)ht(\zeta)\right),$$
where $\deg(W)$ is the degree of the $\Bbb{Q}-$irreducible component $W$
containing $\zeta$. 
\end{theorem}

\noindent The key idea for the proof of this Theorem 
is the use of the $L^3$ (or $LLL$) reduction algorithm.

\spar

\noindent Conversely, as approximate zeros may depend on the height of the actual zero
they approximate, we could be interested in the computation of approximate
zeros for actual zeros with small (bounded) height. 
This is done in the following statement.

\begin{theorem}[From Kronecker's solution to Newton's solution]
\label{theorem-kronecker-newton}
There exists a bounded error probability Turing mahcine $M$ which performs the 
following task~: Given a sequence of polynomial equations 
$F:=(f_1,\ldots,f_n)$ of degree at most $2$ and height at most $h$, 
and given a positive integer number $H\in \mathbb{N}$ in binary encoding, 
the machine $M$ outputs approximate zeros
for the archimedean absolute value $\vert \cdot \vert~: K\longrightarrow \mathbb{R}$
induced on $K$ by the inclusion $i ~: K \hookrightarrow \Bbb{C}$
for all those zeros $\zeta\in V_K(f_1,\ldots,f_n)$, whose logarithmic height
is at most $H$, i.e.
$$ht(\zeta)\leq H.$$
The running time of $M$ is polynomial in the following quantities~:
$$(n\;d\;\delta(F)) + (D\;n\;h\;H),$$
where the notations are the same as in Theorem~\ref{theorem-kronecker-solution}
before. \end{theorem}

\noindent 
A proof of this statement is based  again on an application of the
$L^3$ reduction algorithm.

\section{Application~: Computation of splitting field and Lagrange resolvent}

Combining both Kronecker's and Newton's approach to solving, 
we exhibit an efficient procedure for computing the splitting field and the 
Lagrange resolvent  of an irreducible monic univariate polynomial $f\in
\Bbb{Q}[X]$ of degree $d$.  Let us recall that the splitting field
of $f$ is the minimal number field $K(f)$ containing the field of rational 
numbers $\Bbb{Q}$ and all complex roots of $f$  (i.e. the minimal number field
where $f$ splits completely, also called the normal closure of the equation
$f=0$).  This normal closure $K(f)$ is nothing but the Galois field of $f$ and
it satisfies 
where $Gal_{\Bbb{Q}}(f)$ is the Galois group of the polynomial $f$. 
The splitting field of $f$
can be indentified with an irreducible component of the zero--dimensional
algebra (known as the universal decomposition algebra)
$$A:=\Bbb{Q}[X_1,\ldots, X_d]/(\sigma_0-a_0, \ldots,\sigma_{d-1}-a_{d-1}),$$
where $\sigma_0, \ldots, \sigma_{d-1}$ are the elementary symmetric functions
and $f$ is written as $f(X):=a_0+ a_1 X + \cdots + a_{d-1}X^{d-1} + X^d$. 
Let us also observe that
the Lagrange resolvent is nothing but the Chow (or Cayley) elimination 
polynomial of the zero--dimensional resiude algebra $A/{\frak m}$, where
${\frak m}$ is a well-chosen maximal ideal of $A$. 
Therefore, we can also show the following Theorem as a consequence of the 
comparison between Newton's and Kronecker's approach to solving~:

\begin{theorem}[Splitting Field and Lagrange Resolvent]
\label{theorem-lagrange}
There exists a probabilistic Turing machine, which 
for every given univariate polynomial $f\in\Bbb{Z}[X]$ of degree at most $d$ and 
logarithmic height at most $h$ computes the following items~:
\begin{enumerate}
\item Approximate zeros in $\Bbb{Q}[i]$ of all zeros of $f$,
\item a geometric description of the splitting field $K(f)$ of the polynomial $f$,
\item and the Lagrange resolvent of the equation $f=0$.
\end{enumerate}
The running time of $M$ is polynomial in the following quantities~:
$$\sharp\left(\mbox{Gal}(f)\right)\left(dh\right).$$
\end{theorem}

%
%
%

\newcommand{\etalchar}[1]{$^{#1}$}

\end{document}